\documentclass[12pt]{article}
\usepackage{amsmath}
\usepackage{amssymb}
\usepackage[dvips]{graphicx}
\usepackage{epsfig}
\font\eightmsb=msbm10 scaled 800

\def\eps{\varepsilon}
\font\tencmmib=cmmib10 \skewchar\tencmmib '60
\newfam\cmmibfam
\textfont\cmmibfam=\tencmmib

\font\eightmsb=msbm10 scaled 1100

\def\Bbbb#1{\hbox{\eightmsb#1}}
\def\bbox{\quad\hbox{\vrule \vbox{\hrule \vskip2pt \hbox{\hskip2pt
\vbox{\hsize=1pt}\hskip2pt} \vskip2pt\hrule}\vrule}}
\def\lessim{\ \lower4pt\hbox{$
\buildrel{\displaystyle <}\over\sim$}\ }
\def\gessim{\ \lower4pt\hbox{$\buildrel{\displaystyle >}
\over\sim$}\ }

\def\X{{\cal X}}
\def\A{{\cal A}}
\def\eps{{\varepsilon}}

\def\qed{\hfill\break\rightline{$\bbox$}}
\newcommand{\e}{{\Bbbb E}}
\newcommand{\p}{{\Bbbb P}}
\newcommand{\Reals}{\mathbb{R}}

\newtheorem{proposition}{Proposition}
\newtheorem{lemma}{Lemma}
\newtheorem{theorem}{Theorem}

\makeatletter
\@addtoreset{equation}{section}

\makeatother

\begin{document}

\title{
Deviation inequality for monotonic Boolean functions with application
to a number of $k$-cycles in a random graph.}

\author{
Dmitry Panchenko
\\
{\it Massachusetts Institute of Technology}\\
}
\date{}

\maketitle

\begin{abstract}
Using Talagrand's concentration inequality on the discrete 
cube $\{0,1\}^m$ we show that given a real-valued function $Z(x)$
on $\{0,1\}^m$ that satisfies certain monotonicity conditions
one can control the deviations of $Z(x)$ above its median 
by a local Lipschitz norm of $Z$ at the point $x.$
As one application, we give a simple proof of a nearly optimal
deviation inequality for the number of $k$-cycles
in a random graph. 
\end{abstract}

\vskip 5mm


\section{Introduction and main results.}

In this paper we suggest a new way to use Talagrand's concentration
inequality on the cube to control the deviations of Boolean functions
that satisfy certain monotonicity conditions.
As one application we prove a suboptimal deviation inequality
for the count of $k$-cycles in a random graph.

Let $\X=\{0,1\}$ and define a probability measure 
$\mu$ on $\X$ by $\mu(\{1\})=p, \mu(\{0\})=1-p.$
Consider a product space $\X^m$ with a product probability
measure $\p=\mu^m.$ 
Given a function $Z: \X^m \to\Reals$ and
a point $x=(x_1,\ldots,x_m)\in \X^m$
we define
$$
V_i(x)=Z(x) - Z(x_1,\ldots,x_{i-1},0,x_{i+1},\ldots,x_m),
$$
and 
$$
V(x)=\sum_{i=1}^m V_i^2(x).
$$
Note that $V_i(x)=0$ if $x_i=0.$ 

Let us state the main result of this paper.

\begin{theorem}
If $Z(x)$ and $V_i(x), i\leq m$ 
are non-decreasing in each coordinate 
then for any $a\in \Reals$ and $t>0,$
\begin{equation}
\p\bigl(Z(x)\geq a +\sqrt{V(x) t}\bigr)
\p\bigl(Z(x)\leq a\bigr) \leq e^{-t/2}.
\label{selfnorm}
\end{equation}
\end{theorem}

To understand the statement of Theorem 1, we notice
that a function $(V(x))^{1/2}$ can be interpreted as
a kind of discrete Lipschitz norm of $Z$ locally at the point $x.$
For example,
since $Z(x)$ is defined only on the vertices of $m$-dimensional cube,
if one extends $Z(x)$ linearly from the point
$x$ to its neighbours only along the coordinates where
$x_i=1,$ then $(V(x))^{1/2}$ is the norm of that linear map.
Indeed, if we denote the map by $L: \Reals^m \to \Reals,$ then
$$
L(z)=\sum_{i=1}^m V_i(x)(z_i-x_i) + Z(x),
$$
and $\|L\|=(V(x))^{1/2}.$

The proof of Theorem1 is based on Talagrand's concentration inequality 
on $\X^m.$ In order to give more clear interpretation of (\ref{selfnorm}), 
let us compare it with some typical ways of using Talagrand's inequality. 
One common application is the following.
Given a convex function $f:[0,1]^m\to\Reals$
with a Lipschitz norm 
$$
\|f\|_{L}=\sup_{x,y \in [0,1]^m} \frac{|f(x)-f(y)|}{|x-y|}< \infty,
$$ 
where the supremum is taken over $x\not = y,$
the following inequality holds:
\begin{equation}
\p\bigl(f\geq a + \|f\|_{L}  \sqrt{t}\bigr)
\p\bigl(f\leq a\bigr) \leq e^{-t/2}.
\label{convin}
\end{equation}
If the function $f$ is defined only on the vertices of the cube
$\X^m,$ then it is possible to state a similar result
where one has to use a discrete analog of the Lipschitz norm.
For example, the following deviation inequality holds (see \cite{Bobkov}).
For $i\leq m$ we define
$x^{i}\in \X^m$ such that $x_j^i = x_j$ for $j\not = i$
and $x_i^i \not = x_i,$ and define
$$
\|f\|_{d}=
\sup_{x\in\X^m}\bigl(\sum_{i=1}^m (f(x) - f(x^i))^2\bigr)^{1/2}.
$$
Then
\begin{equation}
\p\bigl(f\geq \e f + \|f\|_{d}  \sqrt{t}\bigr)
\leq  e^{-t/4}.
\label{bob}
\end{equation}
Both inequalities (\ref{convin}) and (\ref{bob})
use global Lipschitz condition
to control the deviation of $f(x).$ Theorem 1 suggests a possibility
of using a local Lipschitz norm $V(x)^{1/2}$ at the point $x,$ 
provided that the monotonicity conditions are satisfied. 
The reason why we compute the Lipschitz norm
only in the direction of decreasing $Z$ is because we control the deviation
of $Z$ {\it above} level $a.$
Theorem 1 is similar in spirit to the ideas in \cite{Kimvu},
\cite{Vu1},\cite{Vu2} (see also references therein),
where the authors describe a way of using
average Lipschitz norm of $Z$ to control its deviations.

One example when the monotonicity conditions are satisfied 
is the following.
Let us consider a set of indices ${\cal M}=\{1,2,\ldots,m\}$
and a set of nonnegative numbers $\alpha_{\cal C}\geq 0$ indexed by
the subsets
${\cal C}\subseteq {\cal M}.$ Consider the function $Z(x)$ defined by
\begin{equation}
Z(x)=\sum_{{\cal C}\subseteq {\cal M}}\alpha_{\cal C}
\prod_{i\in {\cal C}} x_i.
\label{pos}
\end{equation}
In this case the fact that $\alpha_{\cal C}$'s are non-negative
implies that functions $Z(x)$ and $V_i(x), i\leq m$ are non-decreasing
in each coordinate. 
Below we will consider the example of
counting the number of $k$-cycles in a random graph which
can be represented in the form (\ref{pos}) and, thus,
Theorem 1 is applicable.

Consider a standard Erd\"os-R\'enyi
model of a random graph $G(n,p).$ 
Let $V$ be a set of
$n$ vertices, $m={n\choose 2}$ 
and let $E=\{e_1,\ldots,e_m\}$ denote a set
of edges of a complete graph $K_n$ on $n$ vertices.
Given $x=(x_1,\ldots,x_m)\in\X^m,$
the fact that $x_i=1$ or $0$ describes that the edge $e_i$
is present or not present in the graph $G(n,p)$ respectively.  
Let 
$$
C_k=\Bigl\{\{e_{i_1},\ldots,e_{i_k}\} : e_{i_j}\in E, \mbox{ and
$\{ e_{i_1},\ldots,e_{i_k} \}$
form a $k$-cycle} 
\Bigr\}
$$
be a collection of all $k$-cycles, and for $e\in E$ let
$$
C_k(e)=\{c\in C_k : e\in c\}
$$
be a set of all $k$-cycles containing the edge $e.$
We consider the following function on $\X^m$
$$
Z(x)= \sum_{c\in C_k}\prod_{e\in c}x_e,
$$
which is the number of $k$-cycles in a random graph $G(n,p)$.
In this case $V(x)$ can be clearly written as
\begin{equation}
V(x)=\sum_{e\in E} x_e \Bigl(
\sum_{c\in C_k(e)}\prod_{e'\not = e}x_{e'}
\Bigr)^2.
\label{Vee}
\end{equation}
In this case, in order to use Theorem 1 to control the deviation 
of $Z(x)$ above its median $M(Z)$ (or its expectation $\e Z$)
we will proceed by showing how to control
$V(x)$ in terms of $Z(x).$ 

We assume that for some large enough $C(k)>0,$
\begin{equation}
np \geq C(k) \log n.
\label{pee}
\end{equation}
The following theorem holds.
\begin{theorem}
If (\ref{pee}) holds then
there exists a constant 
$C(k)>0$ that depends on $k$ only such that
$$
\p\Bigl(
V(x)\geq C(k)
\bigl( (np)^{k-2}Z(x) + (np)^{2(k-1)} \bigr)
\Bigr) \leq \exp\Bigl(-\frac{(np)^2}{C(k)\log\log np}\Bigr).
$$
\end{theorem}

Theorems 1 and 2 will readily imply the following theorem.

\begin{theorem}
If (\ref{pee}) holds then,

(1) For any $\eps>0$ there exists a constant $C(k,\eps)$
that depends on $k$ and $\eps$ only such that the following holds
$$
\mbox{If }\,\,\, 
\e Z\geq C(k,\eps) 
\,\,\,\mbox{ then }\,\,\, 
M(Z)\leq (1+\eps)\e Z.
$$

(2) There exists a  constant $C(k)>0$ such that
if $\e Z \geq C(k)$ then
\begin{equation}
\p(Z\geq 2\e Z)\leq \exp\Bigl(- \frac{(np)^2}{C(k)\log\log np}\Bigr).
\label{dev}
\end{equation}
\end{theorem}

Recently the authors of \cite{Jan2} proved a more general result
describing the deviations of the count of any subgraph in a random
graph. In the case of $k$-cycles their bound gives
$$
\p(Z\geq 2\e Z)\leq \exp(-(np)^2/C(k)).
$$ 
This shows that
the factor $\log\log np$ in (\ref{dev}) is unnecessary, but
at the moment we don't see how to get rid of it using our approach.
This has nothing to do with Theorem 1, since the factor $\log\log np$
comes directly from Theorem 2 which, probably, can be improved.
In the case of triangles ($k=3$) the bound
$\p(Z\geq 2\e Z)\leq \exp(-(np)^2/C(k))$
was also proved in \cite{Kimvu}.

\section{Proof of Theorem 1.}

Talagrand's concentration inequality on the discrete cube
is the main tool in the proof of Theorem 1.
Let us recall it first.

Given a point $x\in\X^m = \{0,1\}^m$ and a set
${\cal A}\subseteq \X^m,$ let us
denote
$$
U_{\cal A}(x)=\{(s_i)_{i\leq m}\in\{0,1\}^m , 
\exists y\in{\cal A} , s_i=0 \Rightarrow y_i=x_i\}.
$$
The "convex hull" distance between the point 
$x$ and the set $\cal A$ is defined as
$$
f_c({\cal A},x)=\inf\{|s| : s\in\mbox{conv}U_{\cal A}(x)\},
$$
where $|s|$ is the Euclidean norm of $s.$
The concentration inequality of Talagrand (Theorem 4.3.1 in \cite{Ta1})
states the following.

\begin{proposition}
For any $t>0,$
\begin{equation}
\p({\cal A})
\p\Bigl(x\in\X^m : f_{c}^{2}({\cal A},x)\geq t\Bigr)\leq
e^{-t/2}.
\label{T1}
\end{equation}
\end{proposition}

The main feature of this distance is that (Theorem 4.1.2 in \cite{Ta1})
\begin{equation}
\forall (\lambda_i)_{i\leq m}\,\,\,\,\, 
\exists y\in{\cal A}\,\,\,\,\,\,\,\,
\sum_{i=1}^{m}\lambda_i I(y_i\not = x_i)\leq
f_c({\cal A},x)\Bigl(\sum_{i=1}^{m}\lambda_i^2\Bigr)^{1/2}.
\label{T2}
\end{equation}
{\bf Proof of Theorem 1.} For a fixed number $a\in \Reals$ consider a set
$$
\A=\{y\in \X^m : Z(y)\leq a\}.
$$
For a fixed $x\in \X^m$ and an arbitrary $y\in\A,$ 
since $Z(y)\leq a,$ we can write
$Z(x)-a\leq Z(x) - Z(y).$
Consider three sets of indices 
$$
I_1 =\{i: x_i=1,  y_i=0\},\,\,\,\,
I_2 = \{i: x_i = 0, y_i = 1\}, \,\,\,\,
I_3 = \{i: x_i = y_i\}.
$$
Without loss of generality we will assume that $I_1=\{1,\ldots,k\}$
and $I_2 = \{k+1,\ldots, l\}.$
Define a sequence 
$$
z^i = (y_1,\ldots,y_i,x_{i+1},\ldots, x_m), \,\,\,
i=0,\ldots,m.
$$
We have
$$
Z(x) - Z(y)=\sum_{i=1}^{m} (Z(z^{i-1})-Z(z^i))=
\sum_{i=1}^{m} (Z(z^{i-1})-Z(z^i))I(x_i\not = y_i),
$$
since $x_i=y_i$ (i.e. $i>l$) implies that $Z(z^{i-1})-Z(z^i)=0.$
We have 
$$
Z(z^{i-1})-Z(z^i)=0\leq V_i(x)\,\,\,
\mbox{ for }\,\,\,
i=l+1,\ldots,m,
$$
since for this range of indices $z^{i-1} = z^i,$
$$
Z(z^{i-1})-Z(z^i)\leq 0\leq V_i(x)\,\,\,
\mbox{ for }\,\,\,
i=k+1,\ldots,l,
$$
since the function $Z$ is non-decreasing in each coordinate and for 
$i\in I_2,$
$z^{i-1}_i = 0,$ $z_i^i = 1$ and all other coordinates
of $z^{i-1}$ and $z^i$ coincide, and
$$
Z(z^{i-1})-Z(z^i) = V_i(z^{i-1})\leq V_i(x)\,\,\,
\mbox{ for }\,\,\,
i=1,\ldots,k,
$$ 
since for $i\in I_1$
each coordinate of $z^{i-1}$ is smaller than the corresponding
coordinate of $x.$
Thus we proved that for any $y\in {\cal A}$
$$
Z(x) - a \leq \sum_{i=1}^m V_i(x) I(x_i\not =y_i).
$$
By (\ref{T2}) there exists $y\in {\cal A}$ such that
the last expression can  be bounded
$$
\sum_{i=1}^m V_i(x) I(x_i\not =y_i) \leq
f_c({\cal A},x)\Bigl(\sum_{i=1}^m V_i^2(x)\Bigr)^{1/2} =
f_c({\cal A},x)\sqrt{V(x)}
$$
Talagrand's inequality (\ref{T1}) states that
$$
\p(f_c(\A,x)\geq \sqrt{t}) \p(\A)
\leq e^{-t/2},
$$
and, therefore, we finally get
$$
\p(Z(x)\geq a + \sqrt{V(x)t})\p(Z(x)\leq a) \leq e^{-t/2}.
$$
\qed

\section{ Proof of Theorem 2.}

We will denote by $d_v$ the degree of a vertex $v.$
Consider the sequence of sets
\begin{equation}
V_1=\{v: d_v < 16 np\},\,\,\,\,
V_j=\{v : d_v\in [2^{j+2} np, 2^{j+3} np)\},\,\,\,
j\geq 2.
\label{sets}
\end{equation}
We will start by stating several basic facts that will be used in the proof of 
Theorem 2.

\begin{lemma}
If (\ref{pee}) holds then,
\begin{equation} 
\p\Bigl(\exists v: d_v \geq (np)^2\Bigr)\leq e^{-(np)^2/2}.
\label{Basic1}
\end{equation}
\end{lemma}
{\bf Proof.}
For a fixed vertex $v$ its degree
$d_v$ is a sum of $(n-1)$ independent variables with the distribution
$\mu.$ 
Using Bernstein's inequality one can easily check that for $np\geq 4$
$$
\p\Bigl(d_v \geq (np)^2\Bigr)\leq e^{-(np)^2}.
$$
The union bound will produce a factor $n$  
and, therefore, using (\ref{pee}) implies (\ref{Basic1}).
\qed

Thus, with high probability 
we can assume that the degree of each
vertex is bounded by $(np)^2$ and, therefore, 
we can only consider the sets $V_j$ in (\ref{sets})
such that $2^{j+2} \leq np$ and, therefore,
$j\leq \log np.$ 

Next we will bound the cardinality of each $V_j.$

\begin{lemma}
For $C(k)>0$ large enough we have,
\begin{equation}
\p\Bigl(
\exists 2 \leq j\leq \log np\,\,\,\,\,
\mbox{\rm card} V_j \geq \frac{np}{j 2^j \log\log np}
\Bigr) \leq
\exp\Bigl(-\frac{(np)^2}{C(k)\log\log np} \Bigr).
\label{Basic2}
\end{equation}
\end{lemma}
{\bf Proof.}
We copy the proof from \cite{Kimvu} (see equation (12) in section 4.2 there).
For a fixed $j\geq 2,$ assume that
$$
\mbox{card} V_j \geq r = \frac{np}{j 2^j \log\log np}. 
$$
In this case, there exists a set of $r$ vertices each with degree at least
$2^{j+2}np.$ It implies that the number of edges containing exactly
one of these vertices exceeds 
$$
(2^{j+2}np - r)r \geq 2^{j+1} np r,
$$
The probability that such a set of edges exists
is bounded by
\begin{eqnarray*}
& &
{n\choose r}{r(n-r) \choose 2^{j+1}npr} p^{2^{j+1}npr}\leq
\exp\Bigl(
r\log\frac{en}{r} + 2^{j+1}npr \log\frac{er(n-r)}{2^{j+1}npr} +
2^{j+1}npr \log p
\Bigr)
\\
& &
\leq
\exp\Bigl(
r\log\frac{en}{r} + 2^{j+1}npr \log\frac{e}{2^{j+1}}
\Bigr)\leq
\exp\Bigl(
-\frac{(np)^2}{C(k)\log\log np}
\Bigr),
\end{eqnarray*}
where in the last inequality we used the estimate
$$
2^{j+1}npr \log\frac{e}{2^{j+1}}\leq \frac{(np)^2}{2j\log\log np}
\log\frac{1}{2^{j-1}}\leq
-\frac{(np)^2}{C(k)\log\log np},
$$
and the first term $r\log(en/r)$ was negligible compared to the second
term. Taking the union bound over $j\leq \log np,$ we get a factor
$\log np$ in front of the exponent that can be ignored by increasing 
$C(k).$
\qed

Before we will state our next lemma, we need to make one remark about the
proof of Theorem 2. 
Multiplying out the right-hand side of (\ref{Vee})
we observe that $V(x)$ can be written as a sum of terms
$$
x_e \prod_{e'\in c}x_{e'}\prod_{e'\in c'}x_{e'},\,
\mbox{ where } e\in E, c,c'\in C_k(e). 
$$
Each of these terms may appear several times, but, clearly,
the number of appearances will be bounded by $C(k)$ that depends
on $k$ only. Each of these term represents two cycles
that have at least one edge in common. There are many different
isometric configurations of such two cycles but, clearly, the
number of them is bounded by a constant that depends on $k$ only.
Hence, $V(x)$ can be decomposed into the sum of the counts
of such pairs of cycles over different configuration.
With minor modifications it is possible to prove the statement of  
the theorem for each of these configuration. 

We will only look
at the pairs of cycles that have exactly one edge in common. 
Let us denote the number of such pair by $W(x).$
We will identify each pair of cycles with an injection 
$\sigma:\{1,\ldots,2k-2\}\to V(G),$ such that
$(\sigma(1),\sigma(2)\ldots,\sigma(k))$ and 
$(\sigma(k),\sigma(k+1),\ldots,\sigma(1))$ 
are the ordered vertices of these two cycles,
and $\sigma(1)\sigma(k)$ is their only common edge.
Let us denote the set of these injections by $\Sigma_0.$

\begin{lemma}
There is a partition of vertices
$V(G)=F_1\cup\ldots\cup F_{2k-2}$ such that
\begin{equation}
W(x)\leq C(k)\mbox{\rm card}
\{\sigma\in\Sigma_0 :\forall i\,\,\, \sigma(i)\in F_i\}.
\label{Basic3}
\end{equation}
\end{lemma}
{\bf Proof.} See Proposition 1.3 in \cite{Friedgut}.
\qed

Let us denote the set in the statement of Lemma 3 by
\begin{equation}
\Sigma=\{\sigma\in\Sigma_0 :\forall i\,\,\, \sigma(i)\in F_i\}.
\label{sigma}
\end{equation}
{\bf Proof of Theorem 2.}
By Lemma 3 and the discussion preceeding Lemma 3,
all we need to do is to estimate the cardinality of
$\Sigma$ in (\ref{sigma}).
Let us consider the event 
$$
{\cal E}=\Bigl\{
\forall j\geq 2\,\,\,\,
\mbox{\rm card} V_j \leq \frac{np}{j 2^j \log\log np}
\Bigr\}
\bigcup
\bigl\{
\forall v: d_v \leq (np)^2
\bigr\}.
$$
By Lemma 1 and Lemma 2 this event holds with probability at least
$$
1-\exp\Bigl(-\frac{(np)^2}{C(k)\log\log np}\Bigr),
$$
for some $C(k)$ large enough.
From now on we assume that this event occurs. 
For each vertex $v\in F_1$ let us denote
$$
S_l(v)=\{(\sigma(1),\ldots,\sigma(l)) : \sigma\in\Sigma,
\sigma(1)=v\},\,\,\,
l\geq 1.
$$
Let us denote $d_v^{+}=\max(d_v,np).$
We will prove that if ${\cal E}$ occurs than
for $l\geq 2$
\begin{equation}
\mbox{card}S_l(v) \leq C(k)d_v^{+} (np)^{l-2}.
\label{path}
\end{equation}
For $l=2,$ this obviously holds with $C(k)=1.$
We proceed by induction over $l.$ Let us decompose
$$
S_l(v)=\bigcup_{j\geq 1} S_l^j(v), 
$$
where 
$$
S_l^j(v)=\{(\sigma(1),\ldots,\sigma(l))\in S_l(v):
\sigma(l-1)\in V_j\}.
$$
\begin{figure}
\begin{center}
\epsfig{file=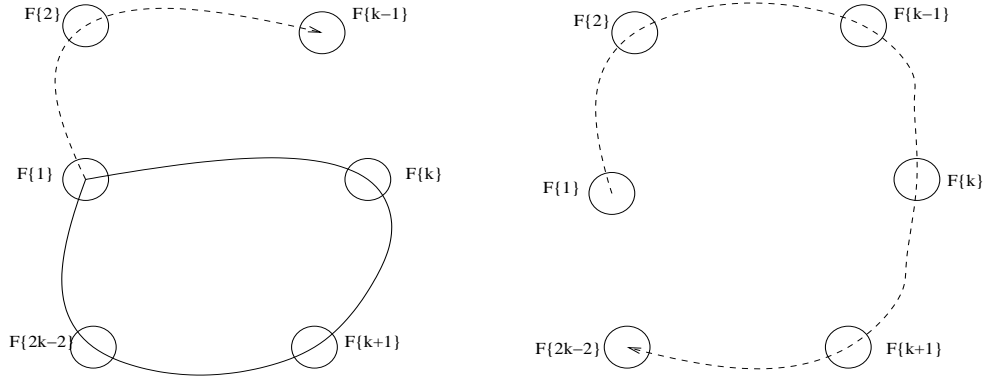,height=5cm,width=13cm}
\caption{
\footnotesize
(1) Counting
$\sigma\in\Sigma_1.$ 
(2) Counting $\sigma\in\Sigma_2.$
}
\label{figure:intervals1}
\end{center}
\end{figure}

To bound the cardinality of $\mbox{card}S_l^1(v)$ 
we use the induction hypothesis and the fact that
the degree of each vertex in the set $V_1$ is bounded by
$16np.$ We get
$$
\mbox{card}S_l^1(v)\leq (16np)\ \mbox{card}{S_{l-1}(v)}\leq
C(k)d_v^{+} (np)^{l-2}.
$$
To bound the cardinality of $S_l^j(v)$ for $j\geq 2$
we notice that for each $(\sigma(1),\ldots,\sigma(l))\in S_l^j(v),$
we have $(\sigma(1),\ldots,\sigma(l-2))\in S_{l-2}(v)$ and 
$\sigma(l-1)\in V_j.$ 
On the event $\cal E$ we can control the cardinality of $V_j$
and, moreover, the degree $d_{\sigma(l-1)}\leq 2^{j+3}np.$  
For $l=3,$ since $\mbox{card}S_1(v)=1,$ we get
$$
\mbox{card}S_3^j(v)\leq 
(\mbox{card}V_j)\ (2^{j+3}np)
\leq 
\frac{np}{2^{j}j\log\log np} 2^{j+3}np
\leq  8 (np)^{2}\frac{1}{j\log\log np} 
$$
and since on the event ${\cal E}$ 
we can assume that $2^{j+2}\leq np,$
which implies that $j\leq \log np,$
we get
$$
\sum_{j=2}^{\log np} \mbox{card}S_3^j(v)\leq
8 (np)^{2}
\sum_{j=2}^{\log np}\frac{1}{j\log\log np}\leq 
C(k) (np)^{2}\leq C(k)d_v^{+}(np),
$$
and this proves the induction step for $l=3.$
The last inequality explains the appearance of the factor
$\log\log np$ in Theorem 3.
Similarly, for $l\geq 4$ we get
\begin{eqnarray*}
&&
\mbox{card}S_l^j(v)\leq 
\mbox{card}S_{l-2}(v)\
(\mbox{card}V_j)\ (2^{j+3}np)
\\
&&
\leq 
C(k)d_v^{+} (np)^{l-4}\frac{np}{2^{j}j\log\log np} 2^{j+3}np
\leq  C(k) d_v^{+} (np)^{l-2}\frac{1}{j\log\log np} 
\end{eqnarray*}
and
$$
\sum_{j=2}^{\log np} \mbox{card}S_l^j(v)\leq
C(k)d_v^{+} (np)^{l-2}
\sum_{j=2}^{\log np}\frac{1}{j\log\log np}\leq 
C(k) d_v^{+} (np)^{l-2}.
$$
This completes the proof of the induction step and (\ref{path}).

To estimate the cardinality of $\Sigma$
we will decompose it into $\Sigma=\Sigma_1\cup\Sigma_2,$
where
$$
\Sigma_1=\{\sigma\in\Sigma : \sigma(1)\in V_1\},\,\,\,\,
\Sigma_2=\{\sigma\in\Sigma : \sigma(1)\in V_j, 2\leq j\}.
$$
We will estimate the cardinality of $\Sigma_1$ and $\Sigma_2$
differently (the idea is illustrated in Figure 1).
First of all, since we can control the cardinality of $V_j$ for 
$j\geq 2,$ we will simply use (\ref{path}) for $l=2k-2$
to compute the number of different paths from $v\in F_1$ to
$F_{2k-2},$ and then add them up. This will give us the bound on
cardinality of $\Sigma_2.$
On the other hand, for $\sigma\in\Sigma_1$ we can represent
it as a cycle on $F_1, F_k,\ldots,F_{2k-2}$ and a path
from $F_1$ to $F_{k-1}.$ In this case, the number of cycles is
bounded by $Z(x),$ and to bound the number of paths we again
use (\ref{path}).

Let first estimate the cardinality of $\Sigma_2.$
First of all by (\ref{path}) for each vertex $v\in V_j \cap F_1$ 
$$
\mbox{card} S_{2k-2}(v)\leq C(k) 2^{j+3}(np) (np)^{2k-4}, 
$$
and, therefore, on the event ${\cal E},$
$$
\mbox{card}\{\sigma\in\Sigma : \sigma(1)\in V_j\} \leq
C(k) 2^{j+3}(np)^{2k-3}
(\mbox{card} V_j)\leq
C(k) 2^{j+3}(np)^{2k-3}\frac{np}{j 2^j\log\log np}.
$$
When we add up these injections over $j\geq 2$ we get
$$
\sum_{j=2}^{\log np} 
 C(k) 2^{j+3}(np)^{2k-3}\frac{np}{j 2^j\log\log np}
\leq
 C(k) (np)^{2k-2}.
$$
This accounts for the second term in the bound of the theorem.

Now consider all injections $\sigma$ such that $\sigma(1)\in V_1$.
Consider the trace of the set of images of the injections 
from $\Sigma$ (in other words, pairs of cycles)
on the set
$$
(V_1\cap F_1) \cup F_k \cup F_{k+1}\cup \ldots \cup F_{2k-2},
$$
i.e.
$$
{\cal P}=
\{(\sigma(1),\sigma(k),\sigma(k+1),\ldots,\sigma(2k-2))
: \sigma\in \Sigma, \sigma(1)\in V_1\}.
$$
First of all, the cardinality of ${\cal P}$ is
bounded by $Z(x),$ since ${\cal P}$ can be identified with
the subset of all cycles in the random graph. 
Moreover, for each $(v_1,v_k,v_{k+1},\ldots,v_{2k-2})\in{\cal P},$
the number of injections $\sigma\in\Sigma$ such that
$\sigma(1)=v_1, \sigma(k) = v_k,\ldots, \sigma(2k-2)=v_{2k-2}$
is bounded by $\mbox{card}S_{k-1}(v_1),$
since all values of the injection are fixed except for
$\sigma(2),\ldots,\sigma(k-1).$
But since $v_1\in V_1$ 
implies that the degree $d_{v_1}\leq 16 np$
and, thus, $d_{v_1}^{+}\leq 16 np,$
we have by (\ref{path})
$$
\mbox{card}S_{k-1}(v_1)\leq C(k)(np)(np)^{k-3}\leq C(k)(np)^{k-2}. 
$$
Therefore, the cardinality of all injections 
such that $\sigma(1)\in V_1$ is bounded by
$$
C(k) (\mbox{card}{\cal P}) (np)^{k-2} \leq C(k) Z(x) (np)^{k-2},
$$
which accounts for the first term in the statement of the theorem.
\qed

\section{Proof of Theorem 3.}
 
Theorem 2 implies that for any $a\in \Reals$ and $t>0,$
\begin{eqnarray*}
&&
\p\Bigl(
Z\leq a +\sqrt{Vt}
\Bigr)
\leq
\p\Bigl(
Z\leq a +\Bigl(C(k)
\bigl((np)^{k-2} Z +(np)^{2(k-1)}\bigr)t\Bigr)^{1/2}
\Bigr)
\\
&&
+
\p\Bigl(
V \geq C(k)((np)^k Z +(np)^{2k})
\Bigr)
\\
&&
\leq
\p\Bigl(
Z\leq a +\Bigl(C(k)
\bigl((np)^{k-2} Z +(np)^{2(k-1)}\bigr)t\Bigr)^{1/2}
\Bigr)
+
\exp\Bigl(
-\frac{(np)^2}{C(k)\log\log np}
\Bigr),
\end{eqnarray*}
which implies that
\begin{eqnarray*}
\p\Bigl(
Z\geq a +\Bigl(C(k) 
\bigl((np)^{k-2} Z +(np)^{2(k-1)}\bigr)t\Bigr)^{1/2}
\Bigr)
&\leq &
\p\Bigl(
Z\geq a +\sqrt{V t}
\Bigr)
\\
& + &
\exp\Bigl(
-\frac{(np)^2}{C(k)\log\log np}
\Bigr).
\end{eqnarray*}
Multiplying both sides by $\p(Z\leq a)$ and using Theorem 1 
with $t=2\eps (np)^2$ we get
\begin{eqnarray}
&&
\p\Bigl(
Z\geq a +\Bigl(C(k) \eps 
\bigl((np)^k Z +(np)^{2k}\bigr)\Bigr)^{1/2}
\Bigr)
\p\Bigl(
Z\leq a
\Bigr)
\nonumber
\\
&&
\leq
\p(Z\leq a)\exp\Bigl(
-\frac{(np)^2}{C(k)\log\log np}
\Bigr)
+ 
e^{-\eps(np)^2}.
\label{final}
\end{eqnarray}
If we take 
\begin{equation}
a=M - \sqrt{C(k) \eps ((np)^k M +(np)^{2k})}
\label{a}
\end{equation}
then, clearly, the following two events are equal
$$
\{Z\geq a +\sqrt{C(k) \eps ((np)^k Z +(np)^{2k})}\} =
\{ Z\geq M\},
$$
and, therefore, (\ref{final}) implies that
for $np\geq C(k)$
for large enough $C(k)>0,$
$$
\p\Bigl(
Z\leq M - \sqrt{C(k) \eps ((np)^k M +(np)^{2k})}
\Bigr) \leq
3e^{-\eps(np)^2}.
$$
Since $\e Z \geq a\p(Z\geq a),$ the choice of $a$
as in (\ref{a}) gives
$$
M - \sqrt{C(k) \eps ((np)^k M +(np)^{2k})}\leq
\e Z \Bigl(1-3e^{-\eps(np)^2}\Bigr)^{-1}.
$$
This, clearly, implies the first statement of Theorem 3.

To prove the second statement we use (\ref{final}) with
$a=M,$ and assume that $np$ is large enough, so that
$M\leq (1+\eps)\e Z.$ Then with probability
at least
$$
1-2e^{-\eps (np)^2} - 
\exp\Bigl(
-\frac{(np)^2}{C(k)\log\log np}
\Bigr)
$$
we have
$$
Z\leq (1+\eps)\e Z +\sqrt{C(k) \eps ((np)^k Z +(np)^{2k})}.
$$
Since $\e Z \sim (np)^k,$
for small enough $\eps$ this implies that
$Z\leq 2\e Z,$ which completes the proof of the second statement
of Theorem  3.
\qed

{\bf Acknowledgment.} We want to thank anonymous referee for
helpful comments, especially, for suggesting
the present formulation of Theorem 1.

\vskip 1mm

\hfill\break
Department of Mathematics\hfill\break
Massachusetts Institute of Technology\hfill\break
77 Massachusetts Avenue, Room 2-181\hfill\break
Cambridge, MA, 02139-4307 \hfill\break
URL: http://www-math.mit.edu/\~{}panchenk\hfill\break
e-mail: panchenk@math.mit.edu\hfill\break

\end{document}